\documentstyle[amssymb,amstex,12pt,a4wide]{article}
\newcommand{\Var}{\operatorname{Var}}
\newtheorem{theorem}{Theorem}
\newtheorem{corollary}{Corollary}
\newtheorem{remark}{Remark}
\newtheorem{lemma}{Lemma}
\newtheorem{Figure}{Figure}

\include{graphicx}

\begin{document}

\title{Limit theory for planar Gilbert tessellations}
\author{Tomasz Schreiber\footnote{Research supported by the Polish Minister of Science and Higher Education grant N N201 385234 (2008-2010)} and Natalia Soja,\\
Faculty of Mathematics \& Computer Science,\\Nicolaus Copernicus University,\\Toru\'n, Poland,\\
\textit{e-mail: tomeks,natas at mat.umk.pl}}
\date{}
\maketitle

\paragraph{Abstract}
 {\it  A Gilbert tessellation arises by letting linear segments (cracks) in ${\mathbb{R}}^2$
 unfold in time with constant speed, starting from a homogeneous Poisson point process of germs in
 randomly chosen directions. Whenever a growing edge hits an already existing one, it stops growing
 in this direction. The resulting process tessellates the plane. The purpose of the present paper is
 to establish law of large numbers, variance asymptotics and a central limit theorem for geometric
 functionals of such tessellations. The main tool applied is the stabilization theory for geometric
 functionals.}

\paragraph{keywords}
{\it Gilbert crack tessellation, stabilizing geometric functionals, central limit theorem, law of large numbers.}                                   


\section{Introduction and main results}
 Let ${\mathcal{X}} \subseteq {\mathbb{R}}^2$ be a finite point set. Each $x \in {\mathcal{X}}$ is
 independently marked
 with a unit length random vector $\hat{\alpha}_x$ making a uniformly distributed angle $\alpha_x \in [0,\pi)$
 with the $x$-axis, which is referred to as the {\it usual marking} in the sequel.
 The collection $\bar{{\mathcal{X}}} = \{ (x,\alpha_x) \}_{x \in {\mathcal{X}}}$
 determines a crack growth process (tessellation) according to the following rules. Initially, at the
 time $t = 0,$ the growth process consists of the points (seeds) in ${\mathcal{X}}$.  Subsequently,
 each point $x \in {\mathcal{X}}$ gives rise to two segments growing linearly at constant unit rate
 in the directions of $\hat{\alpha}_x$ and $-\hat{\alpha}_x$  from $x.$ Thus, prior to any collisions,
 by the time $t> 0$ the seed has developed into the edge with endpoints $x - t\hat{\alpha}_x$ and
 $x + t\hat{\alpha}_x,$ consisting of two segments, say the {\it upper} one $[x,x+t\hat{\alpha}_x]$
 and the {\it lower} one $[x,x-t\hat{\alpha}_x].$ Whenever a growing segment is
 blocked by an existing edge, it stops growing in that direction, without affecting the behaviour
 of the second constituent segment though. Since the possible number of collisions is bounded,
 eventually we obtain a tessellation of the plane. The resulting random tessellation process 
 is variously called the Gilbert model/tessellation, the crack growth process, the
 crack tessellation, and the random crack network, see e.g. \cite{GADK, MM} and
 the references therein.

 Let $G(\bar{{\mathcal{X}}})$ denote the tessellation determined
 by $\bar{{\mathcal{X}}}.$ We shall write $\xi^+(\bar{x},\bar{{\mathcal{X}}}),$ $\; x \in {\mathcal{X}},$
 for the total length covered by the upper segment emanating from $x$ in $G(\bar{{\mathcal{X}}}),$
 and likewise we let $\xi^-(\bar{x},\bar{{\mathcal{X}}})$ stand for the total length of the lower segment
 from $x.$ Note that we use $\bar{x}$ for marked version of $x,$ according to our general
 convention of putting bars over marked objects. For future use we adopt the convention that
 if $\bar{x}$ does not belong to $\bar{{\mathcal{X}}},$
 we extend the definition of $\xi^{+/-}(\bar{x},\bar{{\mathcal{X}}})$ by adding $\bar{x}$ to $\bar{{\mathcal{X}}}$
 and endowing it with a mark drawn according to the usual rules.
 Observe that for some $x$ the values of $\xi^{+/-}$ may be infinite. However, in most
 cases in the sequel ${\mathcal{X}}$ will be a realization of the homogeneous Poisson point process
 ${\mathcal{P}} = {\mathcal{P}}_{\tau}$ of intensity $\tau > 0$ in growing windows of the plane. 
 We shall use the so-called {\it stabilization} property of the functionals $\xi^+$ and $\xi^-,$
 as discussed in detail below,  to show that the construction of $G(\bar{{\mathcal{X}}})$ above
 can be extended to the whole plane yielding a well defined process $G(\bar{{\mathcal{P}}}),$
 where, as usual, $\bar{{\mathcal{P}}}$ stands for a version of ${\mathcal{P}}$ marked as described above. 
 This yields well defined and a.s. finite whole-plane functionals $\xi^+(\cdot,\bar{{\mathcal{P}}})$
 and $\xi^-(\cdot,\bar{{\mathcal{P}}}).$ 

 The conceptually somewhat similar growth process whereby seeds are the realization of a time marked
 Poisson point process in an expanding window of ${\mathbb{R}}^2$ and which subsequently grow radially in all
 directions until meeting another such growing seed, has received considerable attention
 \cite{BY2, Chiu, CQ, CQa, CL, HQR, PY2}, where it has been shown that the number of seeds satisfies
 a law of large numbers and central limit theorem as the window size increases.  In this paper we
 wish to prove analogous limit results for natural functionals (total edge length, sum of power-weighted
 edge lengths, number of cracks with lengths exceeding a given threshold etc.) of the crack tessellation
 process defined by Poisson points in expanding windows of ${\mathbb{R}}^2.$ We will formulate this
 theory in terms of random measures keeping track not only the cumulative values of the afore-mentioned
 functionals but also their spatial profiles. 

 Another interesting class of model bearing conceptual resemblance to Gilbert tessellations
 are the so-called lilypond models which have recently attracted considerable attention
 \cite{CV,DL,HL,HM} and where the entire (rather than just directional) growth is
 blocked upon a collision of a growing object (a ball, a segment etc.) with another one.

 To proceed, consider a function $\phi : [{\mathbb{R}}_+  \cup \{ + \infty \}]^2 \to {\mathbb{R}}$
 with at most polynomial growth, i.e. for some $0 < q < +\infty$ 
 \begin{equation}\label{WZROSTWIELOM}
  \phi(r_1,r_2) =  O\left((r_1+r_2)^q\right).
 \end{equation}
 With $Q_{\lambda} := [0,\sqrt{\lambda}]^2$ standing for the square of area $\lambda$ in ${\mathbb{R}}^2,$
 we consider the {\it empirical measure}
 \begin{equation}\label{MULA}
  \mu^{\phi}_{\lambda} := \sum_{x \in {\mathcal{P}} \cap Q_{\lambda}} \phi\left(\xi^+(\bar{x},\bar{{\mathcal{P}}}),
  \xi^-(\bar{x},\bar{{\mathcal{P}}})\right) \delta_{x/\sqrt{\lambda}}.
 \end{equation}  
 Thus, $\mu^{\phi}_{\lambda}$ is a random (signed) measure on $[0,1]^2$ for all $\lambda > 0.$ The
 large $\lambda$ asymptotics of these measures is the principal object of study in this paper. Recalling
 that $\tau$ stands for the intensity of ${\mathcal{P}} = {\mathcal{P}}_{\tau},$ we define
 \begin{equation}\label{ET}
  E(\tau) := {\mathbb{E}} \phi\left(\xi^{+}(\bar{\bf 0},\bar{{\mathcal{P}}}),\xi^-(\bar{\bf 0},\bar{{\mathcal{P}}})\right).
 \end{equation}  
 The first main result of this paper is the following law of large numbers
 \begin{theorem}\label{PWL}
  For any continuous function $f : [0,1]^2 \to {\mathbb{R}}$ we have
  $$ \lim_{\lambda\to\infty} \frac{1}{\lambda} \int_{[0,1]^2} f d\mu^{\phi}_{\lambda}
   = \tau E(\tau) \int_{[0,1]^2} f(x) dx $$
 in $L^p,\; p \geqslant 1.$
\end{theorem}
 Note that this theorem can be interpreted as stating that $E(\tau)$ is the asymptotic {\it mass per point}
 in $\mu^{\phi}_{\lambda},$ since the expected cardinality of ${\mathcal{P}} \cap Q_{\lambda}$ is
 $\tau \lambda.$   
 To characterize the second order asymptotics of random measures $\mu^{\phi}_{\lambda}$ 
 we consider the pair-correlation functions
 \begin{equation}\label{CXX}
  c_{\phi}[x] := {\mathbb{E}} \phi^2\left(\xi^{+}(x,\bar{{\mathcal{P}}}),\xi^-(x,\bar{{\mathcal{P}}})\right),\;
  x \in {\mathbb{R}}^2
 \end{equation}
 and
 \begin{eqnarray}\label{CXY}
  c_{\phi}[x,y] := {\mathbb{E}} \phi\left(\xi^{+}(x,\bar{{\mathcal{P}}} \cup \{ y \}),\xi^-(x,
  \bar{{\mathcal{P}}} \cup \{ y \})\right) \nonumber \\
 {}\cdot\phi\left(\xi^{+}(y,\bar{{\mathcal{P}}} \cup \{ x \}),
  \xi^-(y,\bar{{\mathcal{P}}} \cup \{ x \})\right) - [E(\tau)]^2.
\end{eqnarray}
In fact, it easily follows by translation invariance that $c_{\phi}[x]$ above does not depend
on $x$ whereas $c_{\phi}[x,y]$ only depends on $y-x.$ 
In terms of these functions we define the asymptotic {\it variance per point} 
\begin{equation}\label{VARASYMPT}
 V(\tau) = c_{\phi}[\bar{\bf 0}] + \tau \int_{{\mathbb{R}}^2} c_{\phi}[\bar{\bf 0},x] dx.  
\end{equation}
Notice that in a special case when function $\phi(\cdot,\cdot)$ is homogeneous of degree $k$ (i.e. for $c \in {\mathbb{R}}$ we have $\phi(cr_1,cr_2)=c^k\phi(r_1,r_2)$) one can simplify (\ref{ET}) and (\ref{VARASYMPT}). Then the following remark is a direct
consequence of standard scaling properties of Gilbert's tessellation construction and those of homogeneous Poisson point processes,
whereby upon multiplying the intensity parameter $\tau$ by some factor $\rho$ we get all lengths in $G(\bar{{\mathcal{P}}})$
re-scaled by factor $\rho^{-1/2}.$ 
\begin{remark}\label{Uwaga1}
For $\phi : [{\mathbb{R}}_+  \cup \{ + \infty \}]^2 \to {\mathbb{R}}$ 
homogeneous of degree $k$ we have
\begin{eqnarray}
E(\tau)=\tau^{-k/2}E(1)\nonumber \\
V(\tau)=\tau^{-k}V(1).
\end{eqnarray}
In other words, $E(\cdot)$ and $V(\cdot)$ are homogeneous of degree $-k/2$ and $-k$, respectively. 
\end{remark}
Our second theorem gives the variance asymptotics for $\mu^{\phi}_{\lambda}.$
\begin{theorem}\label{TWVARASYMPT}
 The integral in (\ref{VARASYMPT}) converges and $V(\tau) > 0$ for all $\tau > 0.$ 
 Moreover, for each continuous $f : [0,1]^2 \to {\mathbb{R}}$ 
 $$ \lim_{\lambda\to\infty} \frac{1}{\lambda} \Var \left[ \int_{[0,1]^2} f d\mu^{\phi}_{\lambda} \right]
     = \tau V(\tau) \int_{[0,1]^2} f^2(x) dx. $$
\end{theorem}
Our final result is the central limit theorem
\begin{theorem}\label{CLT}
 For each continuous $f: [0,1]^2 \to {\mathbb{R}}$ the family of random variables
 $$ \left\{ \frac{1}{\sqrt{\lambda}} \int_{[0,1]^2} f d\mu^{\phi}_{\lambda} \right\}_{\lambda > 0} $$
 converges in law to $ {\cal N}\left(0, \tau V(\tau) \int_{[0,1]^2} f^2(x)dx\right)$ as $\lambda \to \infty.$ 
 Even more, we have
 \begin{equation}\label{SzybkoscZbieznosci}
   \sup_{t \in {\mathbb{R}}}  \left| {\mathbf{P}}\left\{ \frac{\int_{[0,1]^2} f d\mu^{\phi}_{\lambda}}{\sqrt{\Var 
    \left[ \int_{[0,1]^2} f d\mu^{\phi}_{\lambda} \right]}} \leqslant t \right\} - \Phi(t) \right| \leqslant \frac{C (\log \lambda)^6}
     {\sqrt{\lambda}}
  \end{equation}
  for all $\lambda > 1$, where $C$ is a finite constant. 
\end{theorem} 
 
 Principal examples of functional $\phi$ where the above theory applies are
 \begin{enumerate}
  \item $\phi(l_1,l_2) = l_1 + l_2.$ Then the total mass of $\mu_{\lambda}^{\phi}$ coincides with the total length of
           edges emitted in $G(\bar{{\mathcal{P}}})$ by points in ${\mathcal{P}} \cap Q_{\lambda}.$ Clearly, the so-defined
           $\phi$ is homogeneous of order $1$ and thus Remark \ref{Uwaga1} applies.
 \item More generally, $\phi(l_1,l_2) = (l_1+l_2)^{\alpha},\; \alpha \geq 0.$ Again, the total mass of $\mu_{\lambda}^{\phi}$ is
          seen here to be the sum of power-weighted lengths of edges emitted in $G(\bar{{\mathcal{P}}})$ by points in ${\mathcal{P}}
          \cap Q_{\lambda}.$ The so-defined $\phi$ is homogeneous of order $\alpha.$
 \item $\phi(l_1,l_2) = {\bf 1}_{\{ l_1 + l_2 \geq \theta \}},$ where $\theta$ is some fixed threshold parameter. In this set-up,
          the total mass of $\mu_{\lambda}^{\phi}$ is the number of edges in $G(\bar{{\mathcal{P}}})$ emitted from points
          in ${\mathcal{P}} \cap Q_{\lambda}$ and of lengths exceeding threshold $\theta.$ This is not a homogeneous functional.
\end{enumerate}

 The main tool used in our argument below is the concept of {\it stabilization} expressing in geometric
 terms the property of rapid decay of dependencies enjoyed by the functionals considered. The formal
 definition of this notion and the proof that it holds for Gilbert tessellations are given in Section 
 \ref{STABILIZACJA} below. Next, in Section \ref{DOWODY} the proofs of our Theorems \ref{PWL},
 \ref{TWVARASYMPT} and \ref{CLT} are given.

\section{Stabilization property for Gilbert tessellations}\label{STABILIZACJA}
\subsection{Concept of stabilization}
 Consider a generic real-valued translation-in\-var\-iant geometric functional $\xi$ defined on pairs
 $(x,{\mathcal{X}})$ for finite point configurations ${\mathcal{X}} \subset {\mathbb{R}}^2$ and with
 $x \in {\mathcal{X}}.$ For notational convenience we extend this definition for $x \not \in {\mathcal{X}}$
 as well, by putting $\xi(x,{\mathcal{X}}) := \xi(x,{\mathcal{X}} \cup \{ x \})$ then. More generally,
 $\xi$ can also depend on i.i.d. marks attached to points of ${\mathcal{X}},$ in which case
 the marked version of ${\mathcal{X}}$ is denoted by $\bar{{\mathcal{X}}}.$ 

 For an input i.i.d. marked point process $\bar{{\mathcal{P}}}$ on ${\mathbb{R}}^2,$ in this paper
 always taken to be homogeneous Poisson of intensity $\tau,$  we say that the functional 
 $\xi$ {\it stabilizes} at $x \in {\mathbb{R}}^2$ on input $\bar{{\mathcal{P}}}$ iff there exists
 an a.s. finite random variable $R[x,\bar{{\mathcal{P}}}]$ with the property that
 \begin{equation}\label{STABDEF}
  \xi(\bar{x},\bar{{\mathcal{P}}}\cap B(x,R[x,\bar{{\mathcal{P}}}]))=
  \xi(\bar{x},(\bar{{\mathcal{P}}}\cap B(x,R[x,\bar{{\mathcal{P}}}]))\cup \bar{\mathit{A}})
 \end{equation}
 for each finite $\mathit{A} \subset B(x,R[x,\bar{{\mathcal{P}}}])^c,$ with $\bar{\mathit{A}}$ standing
 for its marked version and with $B(x,R)$ denoting ball of radius $R$ centered at $x.$ Note that
 here and henceforth we abuse the notation and refer to intersections of marked point sets with
 domains in the plane \--- these are to be understood as consisting of those marked points whose
 spatial locations fall into the domain considered.  When (\ref{STABDEF})
 holds, we say that $R[x,\bar{{\mathcal{P}}}]$ is a stabilization radius for $\bar{{\mathcal{P}}}$ at $x.$ 
 By translation invariance we see that if $\xi$ stabilizes at one point, it stabilizes at all points of ${\mathbb{R}}^2,$
 in which case we say that $\xi$ {\it stabilizes} on (marked) point process $\bar{{\mathcal{P}}}.$ In addition,
 we say that $\xi$ stabilizes exponentially on input $\bar{{\mathcal{P}}}$ with rate $C > 0$ iff there
 exists a constant $M > 0$ such that 
 \begin{equation}\label{STABWYKL}
  {\mathbf{P}}\{R[x,\bar{{\mathcal{P}}}]>r\}\leqslant Me^{-Cr}
 \end{equation}
 for all $x \in {\mathbb{R}}^2$ and $r> 0.$ Stabilizing functionals are ubiquitous in geometric probability,
 we refer the reader to \cite{BY2,PE2007,PE2007a,PY1,PY2,PY4,PY5,TS09} for further details, where prominent
 examples are discussed 
 including random geometric graphs (nearest neighbor graphs, sphere of influence graphs, Delaunay
 graphs),  random sequential packing and variants thereof, Boolean models and functionals thereof,
 as well as many others.  
 
\subsection{Finite input Gilbert tessellations}
 Let ${\mathcal{X}}\subset{\mathbb{R}}^2$ be a finite point set in the plane. As already mentioned in the
 introduction,  each $x \in {\mathcal{X}}$ is independently marked with a unit length random vector
 $\hat{\alpha}_x = [\cos (\alpha_x),\sin (\alpha_x)]$ making a uniformly distributed angle
 $\alpha_x \in [0,\pi)$ with the $x$-axis and the so marked configuration is denoted by $\bar{{\mathcal{X}}}.$
 In order to formally define the Gilbert tessellation $G(\bar{{\mathcal{X}}})$ as already informally presented above,
 we consider an auxiliary {\it partial tessellation mapping} $G(\bar{{\mathcal{X}}}):{\mathbb{R_{+}}} \to
 {\mathcal{F}}({\mathbb{R}}^2)$ where ${\mathcal{F}}({\mathbb{R}}^2)$ is the space of closed sets in ${\mathbb{R}}^2$
 and where, roughly speaking, $G(\bar{{\mathcal{X}}})(t)$ is to be interpreted as the portion of tessellation
 $G(\bar{{\mathcal{X}}}),$ identified with the set of its edges, constructed by the time $t$ in the course
 of the construction sketched above. 
\begin{Figure}\label{r3}
\begin{center}
\fbox{\includegraphics[width=9cm]{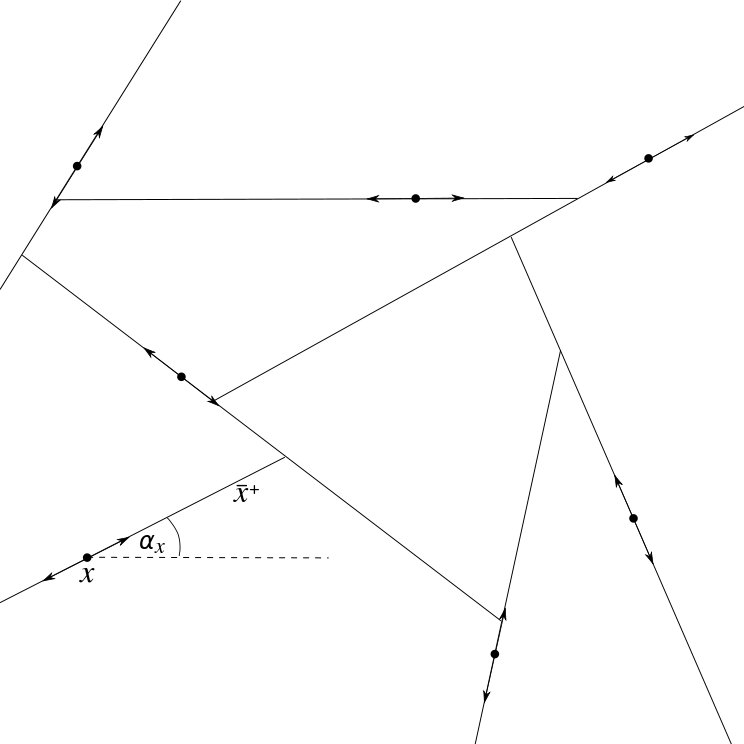}}\\
Finite input Gilbert tessellation.
\end{center}
\end{Figure}
 
  We proceed as follows. For each $\bar{x}=(x,\alpha_x)\in\bar{{\mathcal{X}}}$ at the time moment $0$ the point $x$ 
 emits in directions $\hat{\alpha_x}$ and $-\hat{\alpha_x}$ two segments, referred to as the $\bar{x}^+$- and
 $\bar{x}^-$-branches respectively. Each branch keeps growing with constant rate $1$ in its fixed direction
 until it meets on its way another branch already present, in which case we say it gets {\it blocked}, and
 it stops growing thereupon.  The moment when this happen is called the {\it collision time}. For $t \geqslant 0$
 by $G(\bar{{\mathcal{X}}})(t)$ we denote the union of all branches as grown by the time $t.$ Note that,
 with ${\mathcal{X}} = \{ x_1,\ldots,x_m \},$  the overall number of collisions admits a trivial bound given
 by the number of all intersection points of the family of straight lines  $\{\{x_j+s\hat{\alpha_j}\;  
 s\in{\mathbb{R}}\};j=1,2,\ldots,m\}$ which is $m(m-1)/2$. Thus, eventually there are no
 more collisions and all growth unfolds linearly. It is clear from the definition that  
 $G(\bar{{\mathcal{X}}})(s)\subset G(\bar{{\mathcal{X}}})(t)$ for $s < t.$ The limit set
 $G(\bar{{\mathcal{X}}})(+\infty)=\bigcup_{t\in\mathbb{R_{+}}}G(\bar{{\mathcal{X}}})(t)$
 is denoted by $G(\bar{{\mathcal{X}})}$ and referred to as the {\it Gilbert tessellation}. 
 Obviously, since the number of collisions is finite, the so-defined $G(\bar{{\mathcal{X}}})$
 is a closed set arising as a finite union of (possibly infinite) linear segments. For
 $\bar{x}\in\bar{{\mathcal{X}}}$ by $\xi^+(\bar{x},\bar{{\mathcal{X}}})$ we denote the length
 of the upper branch $\bar{x}^+$ emanating from $x$ and, likewise, we write 
 $\xi^-(\bar{x},\bar{{\mathcal{X}}})$ for the length of the corresponding lower branch.

 For future reference it is convenient to consider for each $x \in {\mathcal{X}}$ the
 {\it branch history}  functions $\bar{x}^+(\cdot),\; \bar{x}^-(\cdot)$ defined by
 requiring that $\bar{x}^{+/-}(t)$ be the {\it growth tip} of the respective branch
 $\bar{x}^{+/-}$ at the time $t \in {\mathbb{R}}_+.$  Thus, prior to any collision in
 the system, we have just $\bar{x}^{+/-}(t) = x +/- \hat{\alpha}_x t,$ that is to say
 all branches grow linearly with their respective speeds $+/- \hat{\alpha}_x.$ 
 Next, when some $\bar{y}^{+/-},\; y \in {\mathcal{X}}$ gets blocked by some
 other $\bar{x}^{+/-},\; x \in {\mathcal{X}}$ at time $t,$ i.e. $\bar{y}^{+/-}(t) =
 \bar{x}^{+/-}(s)$ for some $s \leqslant t,$ the blocked branch stops growing and
 its growth tip remains immobile ever since. Eventually, after all collisions have
 occured, the branches not yet blocked continue growing linearly to $\infty.$  

\subsection{Stabilization for Gilbert tessellations}
 We are now in a position to argue that the functionals $\xi^+$ and $\xi^-$ arising in Gilbert tessellation
 are exponentially stabilizing on Poisson input ${\mathcal{P}} = {\mathcal{P}}_{\tau}$ with i.i.d. marking
 according to the usual rules. The following is the main theorem of this subsection.
\begin{theorem}\label{t1}
 The functionals $\xi^+$ and $\xi^-$ stabilize exponentially on input $\bar{{\mathcal{P}}}.$
\end{theorem}

 Before proceeding to the proof of Theorem \ref{t1} we formulate some auxiliary lemmas.
\begin{lemma}\label{lm1}
 Let ${\mathcal{X}}$ be a finite point set in ${\mathbb{R}}^2$ and $\bar{{\mathcal{X}}}$ the marked
 version thereof, according to the usual rules. Further, let $y \not\in {\mathcal{X}}.$ Then for any
 $t \geqslant 0$ we have
 $$ G(\bar{{\mathcal{X}}})(t) \triangle G(\bar{{\mathcal{X}}}\cup\{\bar{y}\})(t)\subset B(y,t) $$
 with $\triangle$ standing for the symmetric difference.
\end{lemma}
\paragraph{Proof}
 For a point set ${\mathcal{Y}} \subset {\mathbb{R}}^2$ and $x \in {\mathcal{Y}}$ we will use the notation
 $(\bar{x},\bar{{\mathcal{Y}}})^+$ and $(\bar{x},\bar{{\mathcal{Y}}})^-$ to denote, respectively, the
 upper and lower branch outgrowing from $\bar{x}$ in $G(\bar{{\mathcal{Y}}}).$ Also, we use the standard
 extension of this notation for branch-history functions. Note first that, by the construction of
 $G(\bar{{\mathcal{Y}}})$ and by the triangle inequality
 \begin{equation}\label{OBS1}
  (\bar{x},\bar{{\mathcal{Y}}})^{\varepsilon}(s') \in B(y,s') \Rightarrow \forall_{s>s'} 
  (\bar{x},\bar{{\mathcal{Y}}})^{\varepsilon}(s) \in B(y,s),\; s'\geqslant 0, \varepsilon \in \{ -1, +1 \}. 
 \end{equation}
 This is a formal version of the obvious statement that, regardless of the collisions, each branch
 grows with speed at most one throughout its entire history.  

 Next, write ${\mathcal{X}}'={\mathcal{X}} \cup\{ y \}$ and $\Delta(t)=G(\bar{{\mathcal{X}}})(t)\triangle 
 G(\bar{{\mathcal{X}}}')(t)$ for $t\geqslant 0$. Further, let $t_1<t_2<t_3<\ldots<t_n$ be the joint collection of
 collision times for configurations $\bar{{\mathcal{X}}}$ and $\bar{{\mathcal{X}}}'.$ 

 Choose arbitrary $p \in \Delta(t).$ Then there exist unique ${\mathcal{Y}} = {\mathcal{Y}}(p)
 \in \{ {\mathcal{X}}, {\mathcal{X}}'\}$ and $x \in {\mathcal{Y}}$ as well as
 $\varepsilon \in \{+,-\}$ with the property that $p=(\bar{x},\bar{{\mathcal{Y}}})^{\varepsilon}(u)$
 for some $u \leqslant t.$ We also write ${\mathcal{Y}}'$ for the second element of $\{ {\mathcal{X}},
{\mathcal{X}}' \},$ i.e. $\{ {\mathcal{Y}}, {\mathcal{Y}}' \} = \{ {\mathcal{X}}, {\mathcal{X}}' \}.$
With this notation, there is a unique $i = i(p)$ with $t_i$ marking the collision
time in ${\mathcal{Y}}'$ where the branch $(\bar{x},{\mathcal{Y}}')^{\varepsilon}$ gets blocked
in $G(\bar{{\mathcal{Y}}}'),$ clearly $u > t_i$ then  and  for $s<t_i$ we have
 $(\bar{x},\bar{{\mathcal{Y}}})^{\varepsilon}(s)\notin \Delta(t)$.

 We should show that $p\in B(y,t)$. We proceed inductively with respect to $i.$ For $i=0$ we have $x=y$
 and ${\mathcal{Y}} = {\mathcal{X}}'.$
 Since $(\bar{y},\bar{{\mathcal{X}}}')^{\varepsilon}(0)=y\in B(y,0)$, the observation (\ref{OBS1}) implies
 that $p=(\bar{y},\bar{{\mathcal{X}}})^{\varepsilon}(u)\in B(y,u)\subset B(y,t)$. Further, consider the case
 $i > 0$ and assume with no loss of generality that ${\mathcal{Y}}(p)={\mathcal{X}}$, the argument
 in the converse case being fully symmetric. 
 The fact that $p\in G(\bar{{\mathcal{X}}})(t)\triangle G(\bar{{\mathcal{X}}}')(t)$ and that $p=(\bar{x}, 
 \bar{{\mathcal{X}}})^{\varepsilon}(u)$ implies the existence of a point $z\in{\mathcal{X}}'$ such that
 a branch emitted from $z$ does block $\bar{x}^{\varepsilon}$ in $G(\bar{{\mathcal{X}}}')$ (by definition
 necessarily at the time $t_i$) but does not block it in $G(\bar{{\mathcal{X}}})$. In particular, we see that
 $(\bar{z},\bar{{\mathcal{X}}}')^{\delta}(s)=(\bar{x}, \bar{{\mathcal{X}}})^{\varepsilon}(t_i)$ and
 $(\bar{z},\bar{{\mathcal{X}}}')^{\delta}(s')\in \Delta(s')$ for some $\delta,s,s'$ such that $\delta\in\{+,-\}$ 
 and $s'<s\leqslant t_i$. By the inductive hypothesis we get $(\bar{z},\bar{{\mathcal{X}}}')^{\delta}(s')\in B(y,s')$.
 Using again observation (\ref{OBS1}) we conclude thus that $(\bar{x},\bar{{\mathcal{X}}})^{\varepsilon}(t_i)= (\bar{z},\bar{{\mathcal{X}}}')^{\delta}(s)\in B(y,s)$ and hence $p=(\bar{x},\bar{{\mathcal{X}}})^{\varepsilon}(u)\in B(y,u)\subset B(y,t)$. This shows that $p \in B(y,t)$ as required. Since $p$ was chosen arbitrary, this completes the proof of the lemma.
$\Box$

 Our second auxiliary lemma is
\begin{lemma}\label{lm2}
 For arbitrary finite point configuration ${\mathcal{X}} \subset {\mathbb{R}}^2$ and $\bar{x} \in \bar{{\mathcal{X}}}$
 we have
\begin{eqnarray}\label{XISTAB}
\xi^+(\bar{x},\bar{{\mathcal{X}}})=
\xi^+(\bar{x},\bar{{\mathcal{X}}}\cap 
B(x,2\xi^+(\bar{x},\bar{{\mathcal{X}}}))) \nonumber\\
\xi^-(\bar{x},\bar{{\mathcal{X}}})=
\xi^-(\bar{x},\bar{{\mathcal{X}}}\cap
B(x,2\xi^-(\bar{x},\bar{{\mathcal{X}}}))).
\end{eqnarray}
\end{lemma}
\paragraph{Proof}
 We only show the first equality in (\ref{XISTAB}), the proof of the second one being fully analogous.
 Define $A(\bar{{\mathcal{X}}},\bar{x})=\bar{{\mathcal{X}}}\setminus B(x,2\xi^+(\bar{x},\bar{{\mathcal{X}}}))$.  
 Clearly, $A(\bar{{\mathcal{X}}},\bar{x})$ is finite and we will proceed by induction in its cardinality.

 If $|A(\bar{{\mathcal{X}}},\bar{x})|=0$, our claim is trivial. Assume now that $|A(\bar{{\mathcal{X}}},\bar{x})|=n$ 
 for some $n\geqslant 1$ and let $\bar{y}=(y,\alpha_y)\in A(\bar{{\mathcal{X}}},\bar{x})$. Put $t=\xi^+(\bar{x},\bar{{\mathcal{X}}})$ and $\bar{{\mathcal{X}}}'=\bar{{\mathcal{X}}}\backslash\{\bar{y}\}$. 
 Applying Lemma \ref{lm1} we see that $G(\bar{{\mathcal{X}}})(t)\triangle G(\bar{{\mathcal{X}}}')(t)\subset B(y,t)$.  
 We claim that $\xi^+(\bar{x},\bar{{\mathcal{X}}})=\xi^+(\bar{x},\bar{{\mathcal{X}}}')$.
 Assume by contradiction that $\xi^+(\bar{x},\bar{{\mathcal{X}}})\neq\xi^+(\bar{x},\bar{{\mathcal{X}}}')$.
 Then for arbitrarily small $\epsilon>0$ we have $(G(\bar{{\mathcal{X}}})(t)\triangle G(\bar{{\mathcal{X}}}')(t))\cap  B(x,t+\epsilon)\neq\emptyset$. On the other hand, since $\|x-y\|>2t$ as $y \notin B(x,2t),$ for $\varepsilon_0>0$ small enough we get $B(x,t+\epsilon_0)\cap B(y,t)=\emptyset$. Thus, we are led to
$$ 
 \emptyset\neq (G(\bar{{\mathcal{X}}})(t)\triangle G(\bar{{\mathcal{X}}}')(t))\cap B(x,t+\varepsilon_0)\subset B(y,t)
 \cap B(x,t+\varepsilon_0)=\emptyset 
$$ 
which is a contradiction. Consequently, we conclude that 
$t=\xi^+(\bar{x},\bar{{\mathcal{X}}})=\xi^+(\bar{x},\bar{{\mathcal{X}}}')$
as required. Since$|A(\bar{{\mathcal{X}}}',\bar{x})|=n-1$, the inductive hypothesis yields 
$\xi^+(\bar{x},\bar{{\mathcal{X}}}')=\xi^+(\bar{x},\bar{{\mathcal{X}}}'\cap B(x,2\xi^+(\bar{x}, \bar{{\mathcal{X}}}'))=\xi^+(\bar{x},\bar{{\mathcal{X}}}'\cap B(\bar{x},2t))$. Moreover,
$\bar{{\mathcal{X}}}'\cap B(x,2t)=\bar{{\mathcal{X}}}\cap B(x,2t)$. Putting these together we obtain
$$ 
\xi^+(\bar{x},\bar{{\mathcal{X}}})=\xi^+(\bar{x},\bar{{\mathcal{X}}}')
=\xi^+(\bar{x},\bar{{\mathcal{X}}}'\cap B(x,2t))=\xi^+(\bar{x},\bar{{\mathcal{X}}}\cap B(x,2t)),
$$ 
which completes the proof.
$\Box$

 In full analogy to Lemma \ref{lm2} we obtain

\begin{lemma}\label{lm3}
 For a finite point configuration ${\mathcal{X}}\subset{\mathbb{R}}^2$ and $x \in {\mathcal{X}}$ we have
$$ 
 \xi^+(\bar{x},\bar{{\mathcal{X}}})=\xi^+(\bar{x},\bar{{\mathcal{X}}}\cup \bar{A_1})
 \;\; \mbox{ and } \;
 \xi^-(\bar{x},\bar{{\mathcal{X}}})=\xi^-(\bar{x},\bar{{\mathcal{X}}}\cup \bar{A_2})
$$
for arbitrary
$A_1\subset B(x,2\xi^+(\bar{x},\bar{{\mathcal{X}}}))^c$, 
$A_2\subset B(x,2\xi^-(\bar{x},\bar{{\mathcal{X}}}))^c$.
\end{lemma}

 Combining Lemmas \ref{lm2} and \ref{lm3} we conclude

\begin{corollary}\label{w1}
 Assume that finite marked configurations $\bar{{\mathcal{X}}}$ and $\bar{{\mathcal{Y}}}$ coincide on
 $B(x,2\xi^+(\bar{x},\bar{{\mathcal{X}}})).$ Then 
 $$
 \xi^+(\bar{x},\bar{{\mathcal{X}}}\cap 
 B(x,2\xi^+(\bar{x},\bar{{\mathcal{X}}})))=
 \xi^+(\bar{x},\bar{{\mathcal{X}}})=
 \xi^+(\bar{x},\bar{{\mathcal{Y}}}). $$
 Analogous relations hold for $\xi^-.$
\end{corollary}

We are now ready to proceed with the proof of Theorem \ref{t1}.

\paragraph{Proof of Theorem {\ref{t1}}}  
 We are going to show that the functional $\xi^+$ stabilizes exponentially on input process $\bar{{\mathcal{P}}}.$
 The corresponding statement for $\xi^-$ follows in full analogy. Consider auxiliary random variables
 $\xi_{\varrho}^+,\; \varrho >0$ given by
 $$
 \xi_{\varrho}^+=\xi^+(\bar{x},\bar{{\mathcal{P}}} \cap B(x,\varrho))
 $$
 which is clearly well defined in view of the a.s. finiteness of $\bar{{\mathcal{P}}} \cap B(x,\varrho).$
 We claim that there exist constants $M,C>0$ such that for $\varrho\geqslant t\geqslant 0$
\begin{equation}\label{ZANIK1}
 {\mathbf{P}}(\xi_{\varrho}^+>t)\leqslant Me^{-Ct}.
\end{equation}
\begin{Figure}\label{r2}
\fbox{\includegraphics[width=12cm]{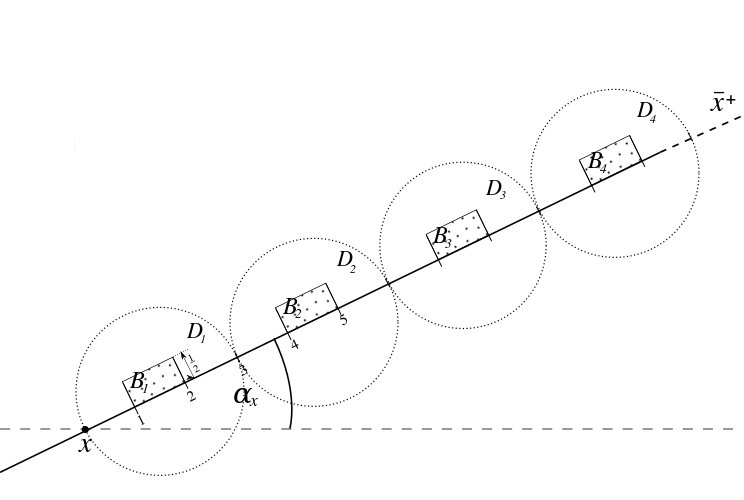}}
\end{Figure}
Indeed, let $\varrho\geqslant 0$. Consider the branch $\bar{x}^+ := (\bar{x},\bar{{\mathcal{P}}} \cap B(x,\varrho))^+$
and planar regions $B_i$ and $D_i,\; i \geqslant 1$ along the branch as represented in figure \ref{r2}.
Say that the event ${\cal E}_i$ occurs iff
\begin{itemize}
 \item the region $B_i$ contains exactly one point $y$ of ${\mathcal{P}}$ and the angular mark $\alpha_y$
          lies within $(\alpha_x + \pi/2 - \epsilon,\alpha_x + \pi/2 + \epsilon),$
 \item and there  are no further points of ${\mathcal{P}}$ falling into $D_i,$
\end{itemize}
where $\epsilon$ is chosen small enough so as to ensure that with probability one on ${\cal E}_i$
the branch $\bar{x}^+$ does not extend past $B_i,$ either getting blocked in $B_i$ or in an
earlier stage of its growth, for instance $\epsilon = 0.01$ will do. Let $p$ stand for the
common positive value of ${\mathbf{P}}({\cal E}_i),\; i \geqslant 0.$ By standard properties of
Poisson point process the events ${\cal E}_i$ are collectively independent. We conclude that,
for ${\mathbb{N}} \ni n \leqslant \varrho / 3$
$$ {\mathbf{P}}(\xi_{\varrho}^+ \geqslant 3 n) \leqslant {\mathbf{P}}\left(\bigcap_{i=1}^n {\cal E}_i^c\right) = (1-p)^n $$
which decays exponentially whence the desired relation (\ref{ZANIK1}) follows. 

Our next step is to define a random variable $R^+=R^+[\bar{x},{\mathcal{P}},\mu]$ and to show it is a 
stabilization radius for $\xi^+$ at $x$ for input process $\bar{{\mathcal{P}}}.$ We shall also establish
exponential decay of tails of $R^+.$ For $\varrho>0$ we put $R^+_{\varrho}=2\xi_{\varrho}^+.$
Further, we set $\hat{\varrho}=\inf \{m\in{\mathbb{N}}\ |\ R^+_m\leqslant m\}.$ Since
${\mathbf{P}}(\bigcap_{m \in {\mathbb{N}}} \{ R^+_m \geqslant m \} ) \leqslant \inf_{m \in {\mathbb{N}}}
 {\mathbf{P}}(R^+_m \geqslant m)$ which is $0$ by (\ref{ZANIK1}), we readily conclude that so defined
$\hat{\varrho}$ is a.s. finite. Take
\begin{equation}\label{PromStab}
 R^+ := R^+_{\hat{\varrho}}.
\end{equation}
Then, using that by definition $R^+ \leqslant \hat{\varrho},$  for any finite $A \subset B(x,R^{+})^c$
we get a.s. by Lemma \ref{lm3} and Corollary \ref{w1}
$$ \xi^+(\bar{x},(\bar{{\mathcal{P}}} \cap B(x,R^+)) \cup A) = \xi^+(\bar{x},\bar{{\mathcal{P}}}
    \cap B(x,\hat{\varrho}) \cap B(x,2\xi^+(\bar{x},\bar{{\mathcal{P}}} \cap B(x,\hat{\varrho}))) \cup A) = $$
$$ =  \xi^+(\bar{x},\bar{{\mathcal{P}}}
    \cap B(x,\hat{\varrho}) \cap B(x,2\xi^+(\bar{x},\bar{{\mathcal{P}}} \cap B(x,\hat{\varrho})))) = 
   \xi^+(\bar{x},(\bar{{\mathcal{P}}} \cap B(x,R^+))). $$
Thus, $R^+$ is a stabilization radius for $\xi^+$ on $\bar{{\mathcal{P}}}$ as required.
Further, taking into account that $R^+_k = R^+$ for all $k \geqslant \hat{\varrho}$ by Corollary \ref{w1}, 
 we have for $m \in {\mathbb{N}}$
\begin{eqnarray}\label{ZANIK2}
  {\mathbf{P}}(R^+ \geqslant m) 
   = {\mathbf{P}}(\lim_{k\to\infty} R^+_k \geqslant m) 
   = \lim_{k \to\infty} {\mathbf{P}}(R^+_k \geqslant m) = \nonumber \\
   = \lim_{k\to\infty} {\mathbf{P}}(\xi^+_k \geqslant m/2)
    \leqslant M e^{-Cm/2} 
\end{eqnarray}
whence the desired exponential stabilization follows.
$\Box$

 Using the just proved stabilization property of $\xi^+$ and $\xi^-$ we can now define
 \begin{equation}\label{XINIESK} 
  \xi^+(\bar{x},\bar{{\mathcal{P}}}) = \xi^+(\bar{x},\bar{{\mathcal{P}}} \cap B(x,R^+)) = 
  \lim_{\varrho\to\infty} \xi^+(\bar{x},\bar{{\mathcal{P}}} \cap B(x,\varrho)) = R^+/2
 \end{equation}
 and likewise for $\xi^-.$ Clearly, the knowledge of these {\it infinite volume} functionals
 allows us to define the whole-plane Gilbert tessellation $G(\bar{{\mathcal{P}}}).$ 

\section{Completing proofs}\label{DOWODY}
 Theorems \ref{PWL},\ref{TWVARASYMPT} and \ref{CLT} are now an easy consequence of
 the exponential stabilization Theorem \ref{t1}. Indeed, observe first that, by (\ref{WZROSTWIELOM}),
 (\ref{XINIESK}) and (\ref{ZANIK2}) the geometric functional 
 $$ \xi(\bar{x},\bar{{\mathcal{X}}}) := \phi(\xi^+(\bar{x},\bar{{\mathcal{X}}}),\xi^-(\bar{x},\bar{{\mathcal{X}}})) $$
 satisfies the p-th bounded moment condition \cite[(4.6)]{TS09} for all $p>0.$ Hence, Theorem \ref{PWL}
 follows by Theorem 4.1 in \cite{TS09}. Further, Theorem \ref{TWVARASYMPT} follows by Theorem
 4.2 in \cite{TS09}. Finally, Theorem \ref{CLT} follows by Theorem 4.3 in \cite{TS09} and Theorem 2.2
 and Lemma 4.4 in \cite{PE2007}.
 
\paragraph{Acknowledgements}
 Tomasz Schreiber acknowledges support from the Polish Minister of Science and Higher Education grant N N201 385234 (2008-2010).
 He also wishes to express his gratitude to J.E. Yukich for helpful and inspiring  discussions.

\end{document}